\documentclass{amsart}

\usepackage{amsmath} 
\usepackage{amsthm}  
\usepackage{amsmath}
\usepackage{amssymb}
\usepackage{amsfonts}
\usepackage{fancyhdr}
\usepackage{hyperref}
\usepackage{subcaption}

\usepackage{bbm}
\usepackage{graphicx}
\usepackage{float}

\newtheorem{lemma}{Lemma}

\newtheorem{theorem}{Theorem}
\newtheorem*{theorem*}{Theorem}

\newtheorem{corollary}{Corollary}

\newtheorem{definition}{Definition}

\newtheorem{proposition}{Proposition}

\newtheoremstyle{nonum}{}{}{\itshape}{}{\bfseries}{.}{.5em}{}
\theoremstyle{nonum}
\newtheorem*{IndexTerms}{Index Terms}

\numberwithin{equation}{section}

\usepackage{graphicx} 

\title{A simplified directional KeRF algorithm}

\author{Iakovidis Isidoros}
\address{Dipartimento di Matematica, Universit\`a di Bologna, 40126, Bologna, Italy}
\email{isidoros.iakovidis2@unibo.it}

\author{Nicola Arcozzi}
\address{Dipartimento di Matematica, Universit\`a di Bologna, 40126, Bologna, Italy}
\email{nicola.arcozzi@unibo.it}

\thanks{The first author was funded by the PhD scholarship PON-DOT1303154-3, and the second author was supported by the Hellenic Foundation for Research and Innovation (H.F.R.I.) under
the ‘2nd Call for H.F.R.I. Research Projects to support Faculty Members \& Researchers (Project Number: 4662).}

\date{March 2025}

\begin{document}

\maketitle

\begin{abstract}
    Random forest methods belong to the class of non-parametric machine learning algorithms. They were first introduced in 2001 by Breiman and they perform with accuracy in high dimensional settings. In this article, we consider a simplified kernel-based random forest algorithm called simplified directional KeRF (Kernel Random Forest). We establish the asymptotic equivalence between simplified directional KeRF and centered KeRF, with additional numerical experiments supporting our theoretical results. Finally, we derive an improved rate of convergence of the centered KeRF in the interpolation regime.

\end{abstract}

\begin{IndexTerms}
Random forests, kernel methods, consistency, non-parametric analysis, randomization.
\end{IndexTerms}

\section{Introduction} 
Random forest algorithms are ensemble machine learning methods for regression and classification tasks. The class of random forest algorithms performs with high accuracy in high-dimensional data, avoiding in general overfitting. It has become one of the most important practical data analysis tools in various active research fields, mainly because of its high performance even in sparse data sets.

Decision trees and random forests have been used in many active research fields in modern sciences and belong to the family of non-parametric machine learning algorithms. In particular, they have been used in economics \cite{Econ}, biology \cite{Biol}, linguistics \cite{Ling}, bio-informatics \cite{Bioinf}, and genetics \cite{CHEN} to name a few. In practice, random forest-type algorithms with many variations are commonly used, since usually, few parameters need to be tuned for training the model.

The random forest algorithm was originally introduced by Breiman in a seminar paper in \cite{BR}. Following the previous work of Ho \cite{Ho}, where it was proposed a method following stochastic modeling for creating multiple tree classifiers, and Amit and Geman \cite{Gem}, who introduced a new approach to shape recognition based on tree classifiers as well, Breiman introduced the CART (Classification and Regression Trees) method. Random forests are widely used in many various research domains since only a few parameters need to be tuned \cite{Li}, \cite{Gen} but still their theoretical properties are under active research \cite{Lin}, \cite{Dev}. For a review of the theoretical and practical aspects of the theory up to 2016, the authors refer to  \cite{BS}.

An important distinction within the class of random forest algorithms lies in the construction of each individual tree. The original random forest algorithm proposed by Breiman \cite{BR} uses the data set to define the estimator through the CART method. Significant research also addresses the consistency of data-dependent algorithms \cite{Ment}, \cite{Wag}, \cite{CONS}. On the contrary, another class of algorithms are those designed independently of the data set the so-called { \it purely random forests}  \cite{Lug}, \cite{AL.SUR2}, \cite{Den}. In 2012, Biau thoroughly examined a specific model proposed by Breiman belonging to the class of { \it purely random forests}, the {\it centered random forest} \cite{Bia}. In this in-depth analysis, Biau proved the consistency of the algorithm and the dependence of the convergence rate only on the strong variables of the feature space and not on the noisy ones.

 Of particular interest to researchers is the rate of convergence of the various random forest algorithms. In particular, in \cite{Bia}, Biau provided an upper bound on the rate of the consistency of the {\it centered random forest}  over the class of the Lipschitz functions of $ \mathcal{O}\left( n^{-\frac{1}{d \frac{4}{3}\log{2}+1}}\right),$ where $n$ is the size of the sample data and $d$ the dimension of the feature space. Moreover, Klusowski in \cite{Sharp_K} provided an improved rate of convergence to $ \mathcal{O}\left( (n\log^{\frac{d-1}{2}} {n})^{ -( \frac{1+\delta}{d\log{2}+1})}  \right),$  where $\delta$ is a positive constant that depends on the dimension of the feature space $d$ and converges to zero as $d$ approaches infinity. The aforementioned rate of convergence is sharp although it fails to reach the minimax rate of convergence in the class of the Lipschitz functions \cite{minimax} $ \mathcal{O}\left(n^{\frac{-2}{d+2}} \right).$

An important direction in manipulating machine learning algorithms is through kernel theory. Already Breiman in \cite{KerB} introduced the connection between kernel methods and random forests where he proved the equivalence between the estimator of a purely random forest and a kernel-type algorithm. In 2006 Geurts et al. in \cite{Geu} formalized the idea by defining the extremely randomized trees method that has a kernel characterization. 

In 2016 Scornet introduced the Kernel Random Forest method (KeRF) \cite{S}. By slightly modifying the definition of random forests, Scornet suggested a new estimator that has a kernel representation. The new KeRF estimators are not of the form of the well-studied Nadaray-Watson kernels with bandwidth and a more careful analysis is necessary. In \cite{S}, the expressions for the centered and uniform KeRF algorithms were explicitly defined, along with upper bounds on their convergence rate. Additionally, due to the kernel representation of the algorithm, the exact formula of the infinite-centered KeRF can be computed (for a more detailed description, see \ref{Kernel}).

Statistical experience shows that high-complex models, where the estimator functions tend to interpolate the data set, easily overfit. Due to the overfitting of the model, statistical intuition suggests that the model tends to perform poorly in new, unseen data, and therefore it fails to generalize appropriately. Kernel interpolation estimators, on the other hand, have been observed to be a good balance between complexity and lack of overfitting \cite{kernel1},\cite{kernel2},\cite{kernel3}. In the work of Belkin et al. \cite{kernel4} non-asymptotic rates with data interpolation were first proven, and recently \cite{rakh} Belkin et al. proved optimal rates of convergence for kernel interpolating estimators. More recently, Wang and Scott \cite{kernel5} provided consistency results for kernel-based methods on Riemannian manifolds.

Arnould et al. \cite{Interp} considered several random forest type algorithms studying the trade off between interpolating estimators and consistency of the algorithms. In particular, they provided rates of convergence in the interpolation regime for the centered KeRF.

\subsection{Contributions of the article}

We introduce a variation of the centered random forest algorithm called the simplified directional algorithm. Our objective is to define a partition of the feature space independently of the data set which simplifies the centered method. We derive the algorithm's corresponding kernel representation and prove that, asymptotically, as the number of trees tends to infinity, the centered KeRF and the simplified directional KeRF coincide. In addition, we support our theoretical results with various experiments, comparing the $L_2$-error and the variance of the centered KeRF and the simplified directional algorithm for different numbers of trees. 

Moreover, we derive improved rates of convergence of the centered KeRF in the mean interpolation regime (see \ref{not}
and \ref{Interpolat} for notation and definitions).

\begin{theorem*}
$   \textbf{Y}=m(\textbf{X}) +\epsilon$ where $\epsilon$ is a zero mean Gaussian noise with finite variance $\sigma $ independent of $\textbf{X} $. Assuming also that $\textbf{X} $ is uniformly distributed in $ [0,1]^d$ and $m$ is a Lipschitz function. Then there exists  constants $c_1, c_2$  depending on $d, \sigma $ and  $ \lVert m \rVert_{\infty} = \sup_{x \in [0,1]^d} | m(x) |$ such that,

\[  \mathbb{E}(\tilde{m}_{n,\infty}^{cc}(x) - m(x))^2 \leq
  c_{1} \left(1-\frac{1}{2d}\right)^{2\log_{2}{n}} 
  +C_3  \frac{ \log_2({ \log_2 n } )^d }{(\log_2 n)^{ \frac{d-1}{2}  }} \log{ 
 \bigg( \frac{ \log{n}  (\log_2 n)^{ \frac{d-1}{2}  }   }{  \log_2({ \log_2 n } )^d    }   \bigg)}    \]

\end{theorem*}
Here, $m(x)=\mathbb{E}\left[Y|X=x\right]$ is the predicted value of $Y$ for $X=x\in [0,1]^d$, while $\tilde{m}^{cc}_{\infty,n}(x)$ is the estimate for $m$ provided by the kernel version of the centered random forest algorithm.

\section{Notation and preliminaries}\label{not}

A fundamental problem in machine learning is to learn the relationship between input features ($\mathbf{X}$) and an output target variable ($Y$) from a set of training examples. Specifically, given a training sample   $(\mathbf{X}_1, Y_1), (\mathbf{X}_2, Y_2), ..., (\mathbf{X}_n, Y_n),$ of independent random variables, where each $\mathbf{X}_i$ is a $d$-dimensional vector in $[0,1]^d$ and $Y_i$ is a real-valued output, the goal is to estimate the function $m(\mathbf{x}) = \mathbb{E}(Y|\mathbf{X}=\mathbf{x})$. In other words, we want to predict the expected value of $Y$ for any new input vector $\mathbf{x}$, based on the training data.
 In classification problems, $Y$ takes on values from a finite set of class labels.

We typically assume that the training sample $(\mathbf{X}_i, Y_i)$ are independent and identically distributed random variables drawn from the same unknown joint distribution $\mathbb{P}(\mathbf{X}, Y)$.

The learning algorithm must then use just these $n$ samples to construct a function $m_n: [0,1]^d \rightarrow \mathbb{R}$ that is a good approximation to the conditional expectation $m(\mathbf{x}) = \mathbb{E}(Y|\mathbf{X}=\mathbf{x})$ for $x \in [0,1]^d.$

Next, we consider the basic framework of the random forest algorithm, which is an ensemble learning method that combines multiple decision trees to improve predictive performance and reduce overfitting.

A {\it random forest} is a finite collection (average) of $M-$ independent, regression random trees designed with respect to the identical distributed random variables $\Theta_{1},..,\Theta_{M}$ as a random variable $\Theta$.

\begin{definition}
For the $j$-th tree in the forest, the predicted value x will be denoted by 
\[ m_{n,\Theta_{j},\mathcal{D}_{n}}(x)=\sum_{i =1}^{n} \frac{\mathbbm{1}_{X_{i}\in A_{n,\Theta_{j},\mathcal{D}_{n} }( x)} Y_{i}}{N_{n,\Theta_{j},\mathcal{D}_{n}}( x) } .\]

Where $ A_{n,\Theta_{j},\mathcal{D}_{n}}( x) $ is the cell containing $x$ and $N_{n,\Theta_{j},\mathcal{D}_{n}}( x)$ is the number of points that fall into the cell that $x$ belongs to.
\end{definition}

For a fixed input $x \in [0,1]^d$, the predicted value of the tree is simply the empirical mean of the target values $Y_i$ for all training points $X_i$ that belong to the same cell as $x$. This is essentially the tree's "guess" for the target value corresponding to $x$, based on the training data.

\begin{definition}
The finite $M$ forest is 
\[ m_{M,n}(x)= \frac{1}{M}\sum_{j=1}^{M} m_{n,\Theta_{j} ,\mathcal{D}_{n}}(x) . \]
\end{definition}
From a modeling point of view, we let $M\to \infty$ and  consider the infinite forest estimate 
\[m_{\infty,n,\mathcal{D}_{n}}(x  )=\mathbb{E}_{\Theta}( m_{n,\Theta,\mathcal{D}_{n} }(x)) . \]
The convergence holds almost surely by the law of the large numbers conditionally on the data set $\mathcal{D}_{n}.$ (Breiman) \cite{AL.SUR}, (Scornet) \cite[Theorem 3.1]{AL.SUR2}.

\subsection{Kernel Random Forest algorithm}\label{Kernel}

In 2016, Scornet introduced kernel methods in the context of random forests, leading to the development of the kernel-based random forest algorithm known as KeRF (Kernel Random Forest) \cite{S}. This innovative approach provided new insights and performance comparisons with traditional random forest algorithms.

To illustrate the underlying principles of KeRF, we begin by revisiting the structure of the traditional random forest algorithm.
For every $x \in [0,1]^d ,$
\[  m_{M,n}(x)= \frac{1}{M}\sum_{j=1}^{M} \big(\sum_{i=1}^n \frac{\mathbbm{1}_{X_{i}\in A_{n,\Theta_{j},\mathcal{D}_{n} }( x)} Y_{i}}{N_{n,\Theta_{j},\mathcal{D}_{n}}( x) }\big). \]

The weights assigned to each observation $Y_{i}$ can be defined as

\[W_{i,j,n}(x)=\frac{\mathbbm{1}_{X_{i}\in A_{n,\Theta_{j},\mathcal{D}_{n} }( x)}}{N_{n,\Theta_{j},\mathcal{D}_{n}}( x) } .\]

The weight values are highly sensitive to the number of points in each cell. The influence of observations within densely populated cells is less compared to those in sparsely populated leaves. To overcome this issue, KeRF considers all tree cells containing a point simultaneously and is defined for all $x \in [0,1]^{d},$ by,
\[  \Tilde{m}_{M,n,\Theta_{1},\Theta_{2},...,\Theta_{M}}(x )=\frac{1}{\sum_{j=1}^{M}N_{n,\Theta_{j}}(x)} \sum_{j=1}^{M}\sum_{i=1}^{n}Y_{i}\mathbbm{1}_{X_{i}\in A_{n,\Theta_{j}}(x)} . \]
This way, empty cells do not affect the computation of the prediction function of the algorithm.

It is proven in \cite{S}, that this representation has indeed a kernel representation.
\begin{proposition}[Scornet \cite{S}, Proposition 1] 
    For all $x \in [0,1]^{d}$ almost surely, it holds
    \[  \Tilde{m}_{M,n,\Theta_{1},\Theta_{2},...,\Theta_{M}}(x )= \frac{\sum_{i=1}^{n}K_{M,n}(x,X_{i})Y_{i}}{\sum_{i=1}^{n}K_{M,n}(x,X_{i})} , \]
    where 
    \[ K_{M,n}(x,z)= \frac{1}{M}\sum_{i=1}^{M}\mathbbm{1}_{ x \in A_{n,\Theta_{i},\mathcal{D}_{n}}( z  ) }. \]
    is the proximity function of the $M$ forest
\end{proposition}
The {\it infinite random forest} arises as the number of trees tends to infinity.
\begin{definition}
   The {\bf infinite KeRF} is defined as:
   \[ \Tilde{m}_{\infty,n}(x)= \lim_{M\to \infty}\Tilde{m}_{M,n}(x,\Theta_{1},\Theta_{2},...,\Theta_{M} ) .  \]
\end{definition}
Infinite KeRF have still the kernel property,

\begin{proposition}[Scornet \cite{S}, Proposition 2]
 Almost surely for all $x,y \in [0.1]^{d}$
\[\lim_{M\to \infty} K_{M,n}(x,y)=K_{n}(x,y), \]
 where
 \[ K_{n}(x,y)= \mathbb{P}_{\Theta} (x\in A_{n}(y,\Theta)),\]
is the probability that $x$ and $y$ belong to the same cell in the infinite forest.
\end{proposition}

\subsection{The Centered KeRF algorithm}

In this paper, we say that an estimator function $m_{n} $ of $m$ is consistent if the following $ L_{2}-$type of convergence holds,
\[\mathbb{E}( m_{n}(x) - m (x) )^{2} \to 0 ,\]
as $n \to \infty.$

In the centered random forest algorithm, the way the data set $\mathcal{D}_{n}$ is partitioned is independent of the data set itself.

\begin{figure}[H]
    \centering
    \includegraphics[width=0.4\textwidth]{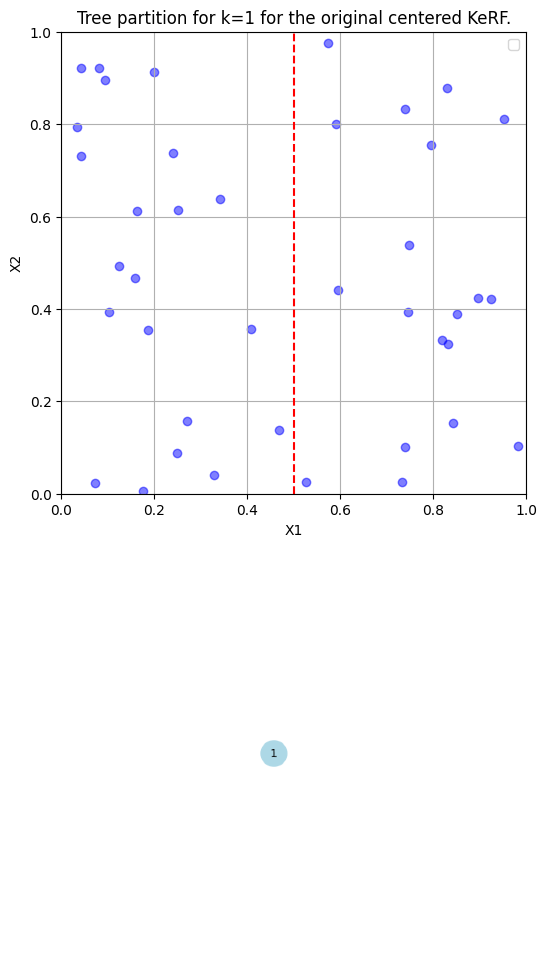}
             \caption{Centered algorithm with tree level  $k=1 $ with the convention that $1$ corresponds to $x_1$ axis and $2$ to the $x_2$ axis.}
    \label{fig:plot_label1}
\end{figure}

\begin{figure}[H]
    \centering
    \includegraphics[width=0.4\textwidth]{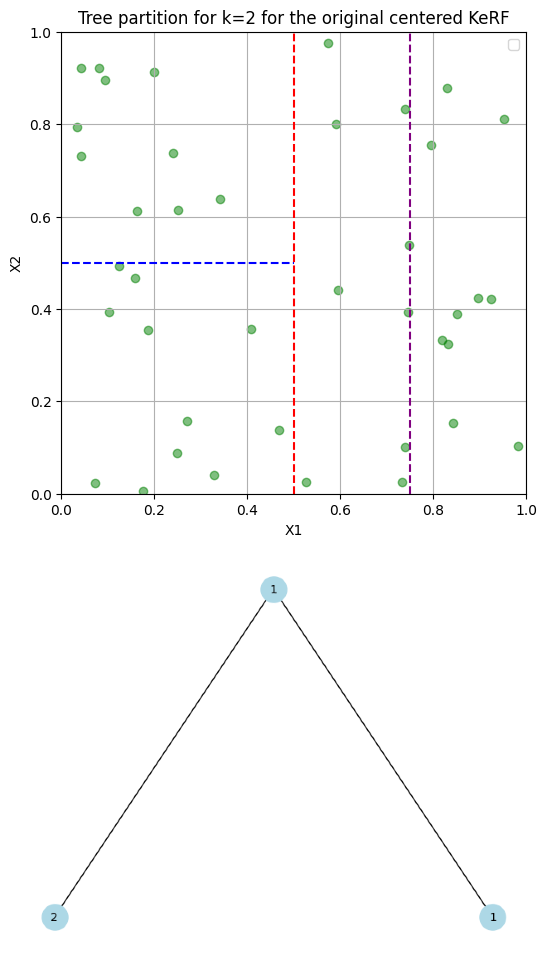}
     \caption{Centered algorithm with tree level  $k=2 $ with the convention that $1$ corresponds to $x_1$ axis and $2$ to the $x_2$ axis.}
    \label{fig:plot_label2}
\end{figure}

The centered forest is designed as follows.
\begin{itemize}
 \item[1)]  Fix $k \in \mathbb{N}.$  
 \item[2)] At each node of each individual tree choose a coordinate uniformly from $\{ 1,2,..d\}.$
 \item[3)] Split the node at the midpoint of the interval of the selected coordinate.
\end{itemize}
Repeat step 2)-3) $k$ times. Finally, we have $2^{k}$ leaves, or cells. A toy example of this iterative process for $k=1,2$ is in Figures \ref{fig:plot_label1},\ref{fig:plot_label2}.
Our estimation at a point $x$ is achieved by averaging the $Y_{i}$ corresponding to the $X_{i}$ in the cell containing $x.$\\

Scornet in \cite{S} introduced the corresponding kernel-based centered random forest providing explicitly the proximity kernel function.

\begin{proposition}
    A centered random forest kernel with $k \in \mathbb{N}$ parameter has the following multinomial expression \cite[Proposition 5]{S},
    \[ K^{Cen}_{k}(x,z) = \sum_{\sum_{j=1}^{d}k_{j}=k } \frac{k!}{k_{1}!...k_{d}!}(\frac{1}{d})^k  \prod_{j=1}^{d}  \mathbbm{1}_{ \left \lceil{2^{k_{j}}x_{j}}\right \rceil =\left \lceil{2^{k_{j}}z_{j}}\right \rceil } . \]
\end{proposition}
Where $K^{Cen}_{k}$ is the Kernel of the corresponding centered random forest. 

\section{Simplified directional KeRF}
In the following section, the simplified directional KeRF algorithm is introduced. We define first the tiling procedure of the hypercube that corresponds to the tree partition and the related forest algorithm. The kernel version of the algorithm is presented and a proof of the asymptotic equivalence of the centered KeRF and the simplified directional KeRF is given.

The simplified directional tree is designed as follows:

\begin{itemize}
 \item[1)]  Fix $k \in \mathbb{N}.$  
 \item[2)] Choose a coordinate uniformly from $\{ 1,2,..d\}.$
 \item[3)] For every node, of each individual tree, split the node at the midpoint of the interval of the preselected coordinate.
\end{itemize}

\begin{figure}[H]
    \centering
    \includegraphics[width=0.4\textwidth]{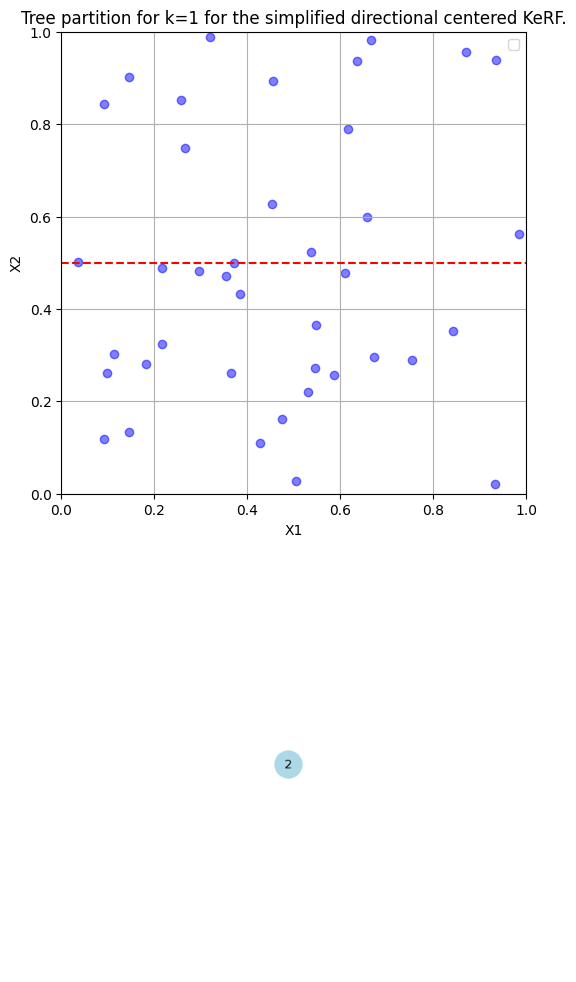}
             \caption{Centered algorithm with tree level  $k=1 $ with the convention that $1$ corresponds to $x_1$ axis and $2$ to the $x_2$ axis.}
    \label{fig:plot_label3}
\end{figure}

\begin{figure}[H]
    \centering
    \includegraphics[width=0.4\textwidth]{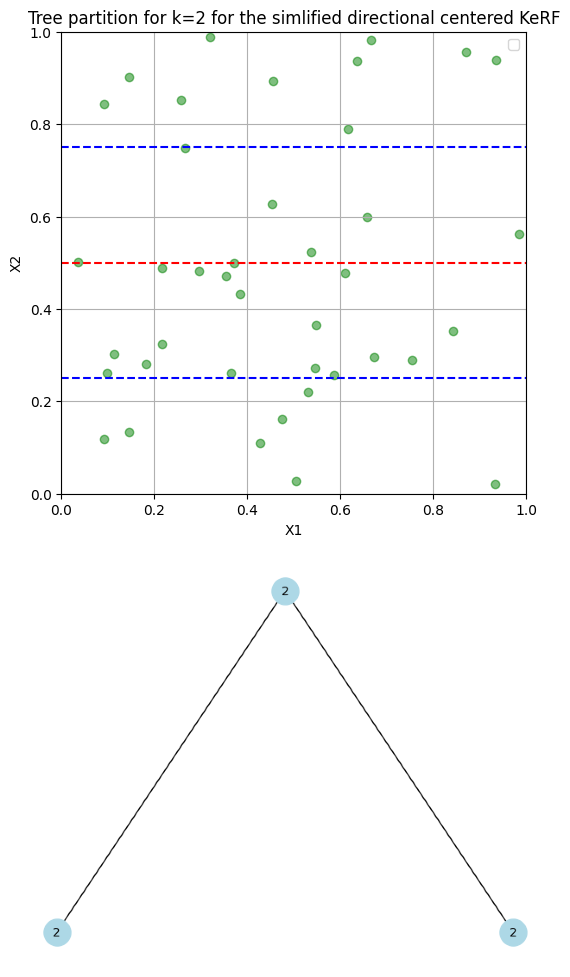}
     \caption{Centered algorithm with tree level  $k=2 $ with the convention that $1$ corresponds to $x_1$ axis and $2$ to the $x_2$ axis.}
    \label{fig:plot_label4}
\end{figure}

Repeat step 2)-3) $k$ times. Finally, we have $2^{k}$ leaves, or cells. A toy example of this iterative process for $k=1,2$ is in Figures \ref{fig:plot_label3},\ref{fig:plot_label4}.

Our estimation at a point $x$ is achieved by averaging the $Y_{i}$ corresponding to the $X_{i}$ in the cell containing $x.$\\

It is clear from the description of the partition of the hypercube for both algorithms, that the latter is indeed a simplification. At each recursive step of the tiling of a tree, in the centered random forest, the choice of the direction of the splitting procedure needs to be taken in every node separately. On the contrary, in the simplified direction random forest, for each recursive step of every tree, there is only a uniform choice for the direction of the splitting.

For simplicity, we define the probability that two points $x,y$ are connected in the $k-$th level of a tree by $p^{sd}_{k} (x,y) $ for the simplified directional algorithm and $p_{k}(x,y)$ respectively for the centered keRF.

\begin{theorem} \label{TH1}
For every $k \in \mathbb{N}$ and every $x,y$ in $[0,1]^d$
\[ p^{sd}_{k} (x,y)=p_{k}(x,y)= K^{Cen}_{k}(x,y)\]
 
\end{theorem}

A simple observation from the theorem above is, that the infinite centered KeRF coincides with the infinite simplified directional KeRF since they have the same infinite kernel representations.
Hence, under specific assumptions, we can compute the rate of convergence of the infinite simplified KeRF. We provide as a simple corollary the speed of convergence over the class of the $L-$ Lipschitz functions and in the chapter \ref{Interpolat} the improved convergence rate of the infinite centered KeRF (and therefore of the infinite simplified KeRF) in the interpolation regime.

\begin{corollary}\label{RATE} (\cite{I-A} Theorem 1.)
Under the following assumptions:
\begin{align*}
    & Y = m(X) + \epsilon  \\
    & X \text{ is uniformly distributed on } [0,1]^{d}  \\
    & \epsilon \sim \mathcal{N}(0,\sigma^2), \quad \sigma < \infty  \\
    & m \text{ belongs to the class of } L\text{-Lipschitz functions,} 
\end{align*}
the rate of convergence of the simplified directional KeRF is  $\mathcal{O} \left( n^{-\bigl(\frac{1}{1+d \log{2}}\bigr)} (\log{n})\right)$

\end{corollary}

\section{Plots and experiments.}

In this section, we conduct numerical simulations and experiments to compare the performance of the finite-centered KeRF algorithm and the finite-simplified directional KeRF algorithm. The evaluation is carried out in terms of the $ L_2$-error and the standard deviation of the error for various target functions $Y$. Specifically, we generated a two-dimensional feature space of size $n=1500$ comprising uniformly distributed points.

The dataset was split into training and testing subsets, with $80\%$  of the data utilized for training both algorithms, while the remaining $20\%$ was reserved for evaluation purposes computing $(  \sum_{ X_{i} \in \; \text{test set} }  (\tilde{m}(X_i) -Y_i)^2 ) $.
To evaluate the performance of the algorithms, we considered several target functions $Y$. For each function, we trained both the finite-centered KeRF and the simplified directional KeRF on the training subset and subsequently evaluated their predictions on the remaining testing subset.

The following target functions with linear, polynomial, and exponential relationships within the feature space were considered to investigate the $L_2$-error with a fixed tree depth value of $k=\log_{2}{n}$ and hence, every leaf has on average $1$ data point. Moreover, the number of trees varies from $ M=1,50,100,200,300,400,500$ and therefore we can empirically confirm that the two algorithms coincide asymptotically.

\begin{enumerate}
\item  $Y=X_1+X_2+\epsilon$, where $\epsilon$ is a random error following a normal distribution $\mathcal{N}(0,\frac{1}{2}).$
\item $Y=X_1^2+X_2^2+\epsilon$, where $\epsilon$ is a random error following a normal distribution $\mathcal{N}(0,\frac{1}{2}).$
\item $Y=2X_1+e^{-X_2^2}.$
\end{enumerate}

All numerical simulations were conducted using the open-source Python software \url{ https://www.python.org/,}, primarily utilizing the numpy library.

\begin{figure}[H]
    \centering
    \begin{subfigure}[b]{0.4\textwidth}
        \centering
        \includegraphics[width=\textwidth]{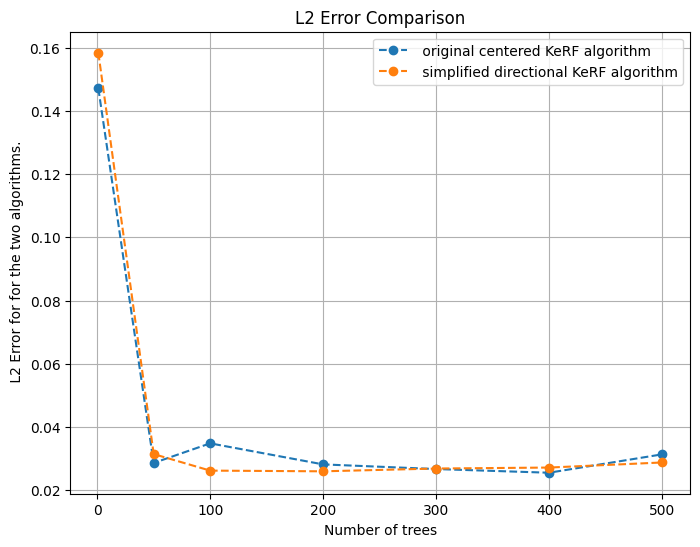}
        \caption{$L-$error for the function $Y=2X_1+e^{-X^{2}_{2}}$.}
        \label{fig:plot_label5}
    \end{subfigure}
    \hfill
    \begin{subfigure}[b]{0.4\textwidth}
        \centering
        \includegraphics[width=\textwidth]{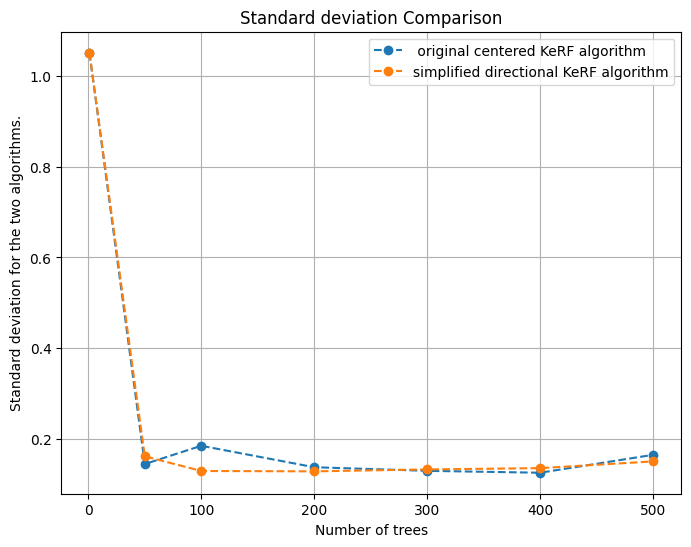}
        \caption{Standard deviation for the $L-$error for the function $Y=2X_1+e^{-X^{2}_{2}}$.}
        \label{fig:plot_label6}
    \end{subfigure}
    \caption{Comparison of $L-$error and standard deviation.}
    \label{fig:combined_label11}
\end{figure}

\begin{figure}[H]
    \centering
    \begin{subfigure}[b]{0.4\textwidth}
        \centering
        \includegraphics[width=\textwidth]{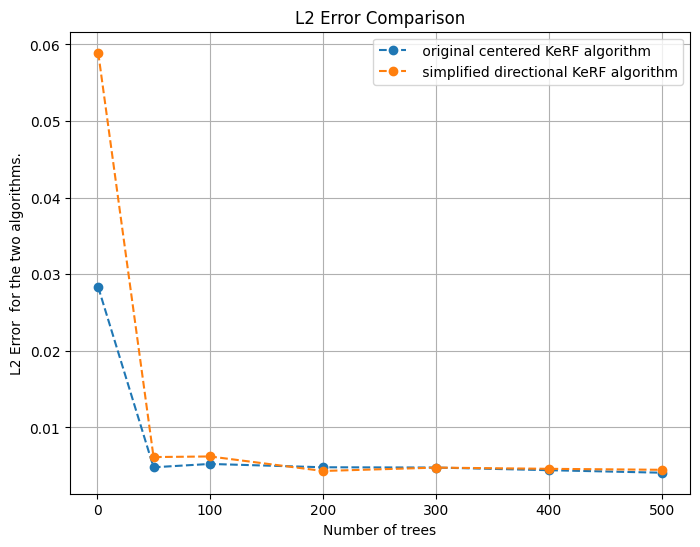}
        \caption{$L-$error for the function $Y=X^{2}_1+ X^{2}_{2} +\epsilon$ where $ \epsilon\sim \mathcal{N}(0,\frac{1}{2})  $.}
        \label{fig:plot_label7}
    \end{subfigure}
    \hfill
    \begin{subfigure}[b]{0.4\textwidth}
        \centering
        \includegraphics[width=\textwidth]{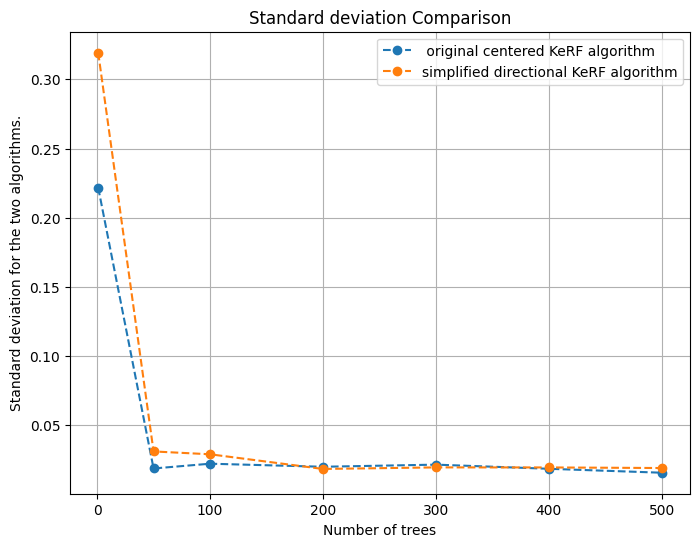}
        \caption{Standard deviation for the $L-$error for the function $Y=X^{2}_1+ X^{2}_{2} +\epsilon$ where $ \epsilon\sim \mathcal{N}(0,\frac{1}{2})  $.}
        \label{fig:plot_label8}
    \end{subfigure}
    \caption{Comparison of $L-$error and standard deviation.}
    \label{fig:combined_label12}
\end{figure}

\begin{figure}[H]
    \centering
    \begin{subfigure}[b]{0.4\textwidth}
        \centering
        \includegraphics[width=\textwidth]{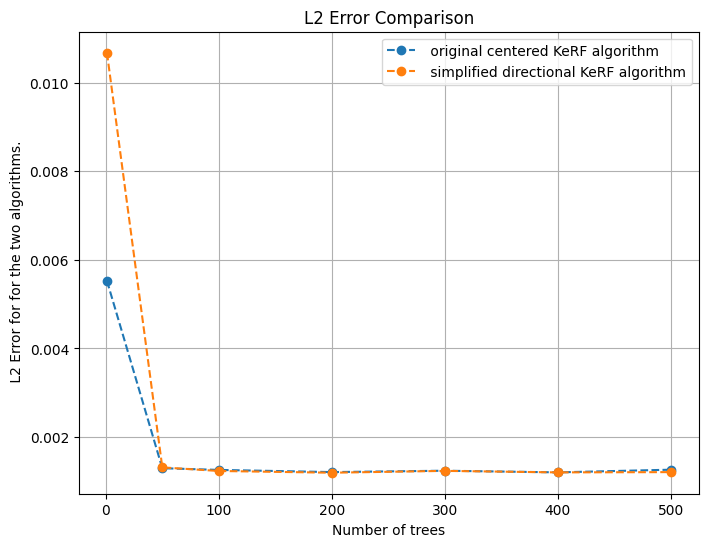}
        \caption{$L-$error for the function $Y=X_1+ X_{2} +\epsilon$ where $ \epsilon\sim \mathcal{N}(0,\frac{1}{2})  $.}
        \label{fig:plot_label9}
    \end{subfigure}
    \hfill
    \begin{subfigure}[b]{0.4\textwidth}
        \centering
        \includegraphics[width=\textwidth]{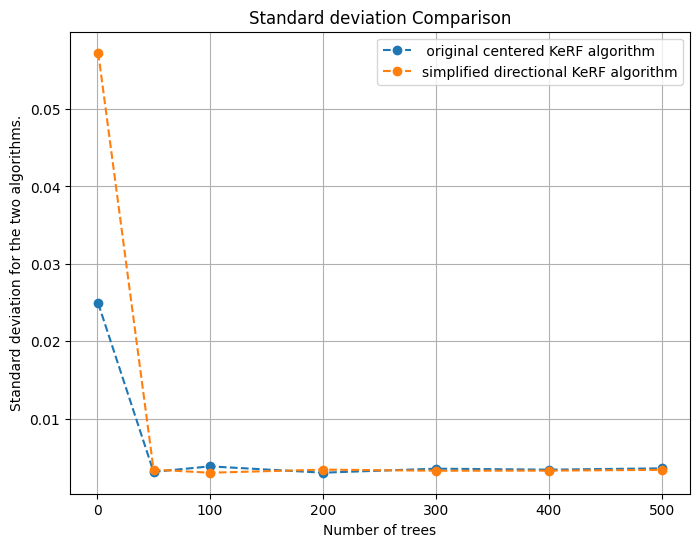}
        \caption{Standard deviation for the $L-$error for the function $Y=X_1+ X_{2} +\epsilon$ where $ \epsilon\sim \mathcal{N}(0,\frac{1}{2})  $.}
        \label{fig:plot_label10}
    \end{subfigure}
    \caption{Comparison of $L-$error and standard deviation.}
    \label{fig:combined_label13}
\end{figure}

As one might expect, all the figures \ref{fig:plot_label5} \ref{fig:plot_label7}, \ref{fig:plot_label9} exhibit similar behavior. For small values of trees, the two algorithms demonstrate slightly different performances, however, as the number of trees increases, consistent with the theorem \ref{TH1}, both algorithms have the same performance in terms of the $L_2$-error consistently with \ref{TH1}. Similarly, the same results hold for the standard deviation of the errors \ref{fig:plot_label6} \ref{fig:plot_label8}, \ref{fig:plot_label10}. Overall, as it is evident from all experiments \ref{fig:combined_label11}, \ref{fig:combined_label12} and \ref{fig:combined_label13}, after $M=100$ trees the centered KeRF and the simplified directional KeRF essentially coincide.

\section{Interpolating random trees}\label{Interpolat}

In this section, we provide an improvement of the rate of convergence of the infinite centered KeRF in the interpolation regime. Since the infinite directional KeRF coincides with the infinite centered KeRF we get simultaneously the rates of convergence for both algorithms. In \cite{Interp} Arnould et al. examined the interpolation of data of several random forest models and their capability to preserve consistency.

Next, we provide the basic definitions on data interpolation and we mention some classical results. Finally, we provide the improved convergence rate.

\begin{definition}
    An estimator $m_n$ interpolates the data set if for every $(X_i, Y_i)$ in the training set we have $m_n(X_i)=Y_i$ almost surely.
\end{definition}

Moreover, since a random forest is an average of $M-$ random trees it is sufficient for the random forest estimator to interpolate the data if every tree estimator interpolates the data set. The tree estimator for a point $x$ is the average of $Y_i$'s for those $X_i$'s belonging to the same cell (or node). Therefore, it is clear that a tree interpolates the data set if a tree is grown until each node contains $X_i$'s with the same values of $Y_i$'s. The regression models we study contain a Gaussian random noise $\epsilon$ and therefore, almost surely, all $Y_i$'s have different values in our data set.

In fact, since $X \text{ is uniformly distributed on } [0,1]^{d}$ the probability that point belongs to one node is $\frac{1}{2^k}$ and the expected number of points per node are $\frac{n}{2^k}.$

\begin{definition}
    A centered random forest estimator satisfies the mean interpolation regime if every tree has at least $n-$ nodes. In other words, if $2^k \geq n .$
\end{definition}
    
The centered random forest fails to preserve consistency in the interpolation regime. This is a result by Arnould et al. \cite{Interp}.
\begin{theorem}[Inconsistency of Centered Random Forest]
If $E[m(X)^2] > 0$ and $k_n \geq \log_2(n)$, then the infinite centered random forest $m^{cc}_{\infty,n}$ is inconsistent.
\end{theorem}

The kernel-based centered algorithm differs from the centered forest only on the way the observations are averaged and not on the way the tiling is performed. Hence, the centered KeRF and the simplified directional KeRF satisfy the mean interpolation regime again if and only if every tree is grown until $2^k \geq n.$
In the same article, Arnould et al. prove that centered kernel random forest are consistent in the interpolation regime with a slow convergence rate, as long as the dimension of the feature space is greater than $5.$ 

Intuitively the reason why the kernel-based centered (or simplified directional) estimator is consistent, although the tree construction of both algorithms is the same, is the way the kernel estimator is computed. By default, the number of empty nodes in each tree partition is the same for a centered random forest and a centered KeRF. The kernel estimator though, does not take into account empty cells, and this is why exactly the consistency is preserved, even with slow convergence rates and deep tree depth. 

Finally, we mention here the theorem of Arnould et al. \cite{Interp} about the convergence rate of the centered KeRF algorithm in the mean interpolation regime and after we state an improved convergence rate that holds even for low dimensional feature spaces.

\begin{theorem}[Consistency of Centered KeRF]
Under the following assumptions:
\begin{align*}
    & Y = m(X) + \epsilon, \\
    & X \text{ is uniformly distributed on } [0,1]^{d}, \\
    & \epsilon \sim \mathcal{N}(0, \sigma^2), \quad \sigma < \infty, \\
    & m \text{ belongs to the class of } L\text{-Lipschitz functions,}
\end{align*}
and assuming furthermore that $k = \lfloor \log_2(n) \rfloor$:

then the rate of convergence is
\[
\mathbb{E}[(\tilde{m}_{n,\infty}^{cc}(x) - m(x))^2] \leq \frac{8L^2d^2}{n^{-2\log_2(1 - 1/d)}} + C_d (\log_2 n)^{-(d-5)/6} (\log_2(\log_2 n))^{d/3},
\]
where $C_d > 0$ is a constant dependent on noise variance.
\end{theorem}

Under the same assumptions for the regression function $m$, in the mean interpolation regime, we provide an improvement in the rate of convergence.

\begin{theorem}
$   \textbf{Y}=m(\textbf{X}) +\epsilon$ where $\epsilon$ is a zero mean Gaussian noise with finite variance $\sigma $ independent of $\textbf{X} $. Assuming also that $\textbf{X} $ is uniformly distributed in $ [0,1]^d$ and $m$ is a Lipschitz function. Then there exists  constants $c_1, c_2$  depending on $d, \sigma $ and  $ \lVert m \rVert_{\infty} = \sup_{x \in [0,1]^d} | m(x) |$ such that,

\[  \mathbb{E}(\tilde{m}_{n,\infty}^{cc}(x) - m(x))^2 \leq
  c_{1} \left(1-\frac{1}{2d}\right)^{2\log_{2}{n}} 
  +C_3  \frac{ \log_2({ \log_2 n } )^d }{(\log_2 n)^{ \frac{d-1}{2}  }} \log{ 
 \bigg( \frac{ \log{n}  (\log_2 n)^{ \frac{d-1}{2}  }   }{  \log_2({ \log_2 n } )^d    }   \bigg)}    \]

\end{theorem}
The above theorem states that the rate of convergence of the infinite centered KeRF and the infinite simplified directional KeRF is faster than the one in \cite{Interp} in the interpolation regime even for the dimension of the feature space $d\geq 2.$ Moreover, interpolation in probability and consistency holds simultaneously but in a relatively slow convergence rate. One can obtain the rate of \ref{RATE} by optimizing the depth parameter.

Lin and Jeon provided a theoretical lower bound for the rate of convergence of deep non-adaptive random forests \cite{Lin} of $\frac{1}{(\log{n})^{d-1}}$ and therefore we do not know yet if our rate of convergence is generally improvable. On the contrary, kernel estimators of the Nadaraya–Watson type (\cite{nad}, \cite{watc}) where the smoothing parameter is highly related with the tree depth parameter have been studied in \cite{rakh} by Belkin and Rakhlin where it was proved that the rate of convergence is the minimax over the class of the Lipschitz functions.

\section{Conclusions}
In this article, we introduce a simplified version of the centered KeRF algorithm, the so-called simplified directional KeRF. We study their kernel-based versions and prove that the two algorithms coincide asymptotically. Moreover, we provide numerical simulations for low-dimensional examples that corroborate our theoretical results and we provide an improved consistency rate in the mean interpolation regime for the centered KeRF.

\section{Proofs of theorems}\label{Proofs}

In this final section, we provide the proofs of all technical theorems and lemmata of the previous sections.

\begin{theorem} \label{TH1}
For every $k \in \mathbb{N}$ and every $x,y$ in $[0,1]^d$
\[ p^{sd}_{k} (x,y)=p_{k}(x,y)= K^{Cen}_{k}(x,y)\]
 \begin{proof}
 For $k=0,1$ the result is trivial. 
 
 We assume that for every $x,y \in [0,1]^d$ and for $l=0,1,...k$, $ p^{sd}_{k} (x,y)=p_{k}(x,y) $ and the proof without loss of generality is provided for $d=2.$ 
 Moreover, let $n_{sd}(k) $ resp $( n(k) ) $ to be the number of different tree expansions of level $k$ for the simplified directional algorithm (resp original centered algorithm), and recursively it is easy to check that
 
 \[n_{sd}(k)=2 n_{sd}(k-1)=...=2^{k} \]
 
 and with the same arguments,
 
 \[ n(k)=2^{k}n(k-1)=...=2^{\frac{k(k+1)}{2}}.\]

 Furthermore, we denote by $A_{x,y}^{k}$ the number of times that the points $x,y$ fall in the same cell in the original centered algorithm and $B_{x,y}^{k} $ for the simplified directional respectively.
 Then, $p^{sd}_{k} (x,y)=\frac{B_{x,y}^{k}}{n_{sd}(k)} $ and $p_{k} (x,y)=\frac{A_{x,y}^{k}}{n(k)} $ and we observe the following cases for the original centered random forest algorithm:

 If $x,y$ are not connected for every possible different tree expansion of level $k$ then \[p_{k} (x,y)=p^{sd}_{k} (x,y)=p_{k+1} (x,y)=p^{sd}_{k+1} (x,y)=0.\] Furthermore, if $x,y$ are connected for some possible different tree expansion of level $k,$ but they are not connected for any tree expansions of level $k+1$ then 

\[ p_{k+1} (x,y)=p^{sd}_{k+1} (x,y)=0.\]
Next, if $x,y$ are connected for some possible different tree expansions of level $k,$ and they are connected only after a horizontal cut and not after a vertical, then 
\begin{align*}
    p_{k+1} (x,y) &= \frac{A^{k+1}_{x,y}}{n(k+1)}\\
    &=\frac{ \frac{1}{2}  2^{k+1}A^{k}_{x,y}  }{2^{\frac{(k+1)(k+2)}{2} }} \quad \text{ since half of the tree expansions are connected}\\
    &=\frac{1}{2}p_{k}\\
    &=\frac{1}{2} p^{sd}_{k} (x,y) \quad \text{by the induction hypothesis}
\end{align*}
 and of course,
 \[ p^{sd}_{k+1} (x,y)=\frac{B^{k+1}_{x,y}}{n_{sd}(k+1)}=\frac{1}{2} \frac{2 B^{k}_{x,y}}{2^{k+1}}=\frac{1}{2}p^{sd}_{k} (x,y).\]
By symmetry, the result holds as well if  $x,y$ are connected for some possible different tree expansion of level $k,$ and they are connected only after a vertical cut.
 Finally, when $x,y$ are connected for some possible different tree expansions of level $k,$ and they are connected as well, after the next cut, in any direction then,

 \begin{align*}
    p_{k+1} (x,y) &= \frac{A^{k+1}_{x,y}}{n(k+1)}\\
    &=\frac{  2^{k+1}A^{k}_{x,y}  }{2^{\frac{(k+1)(k+2)}{2} }} \\
    &=p_{k}\\
    &=p^{sd}_{k} (x,y) \quad \text{by the induction hypothesis}\\  
\end{align*}
and,
\[  p^{sd}_{k+1} (x,y)= \frac{B^{k}_{x,y}}{2^{k}} =p^{sd}_{k} (x,y) ,\] which concludes the proof.

 \end{proof}
\end{theorem}

We assume furthermore that all random variables are real-valued and  $\lvert \lvert X \rvert\rvert_{L_{p}} \colon= (\mathbb{E}\lvert X \rvert^{p} )^{\frac{1}{p}} $ and $ \lvert \lvert X  \rvert \rvert_{\infty} \colon= \inf\{ t \colon P(\lvert X \rvert\leq t )=1 \}  $

We begin with this basic lemma providing tail bounds for centered iid random variables with bounded variance and supremum norm.

\begin{lemma}\label{Bound_lemma1}
  Let $X_{1},..., X_{n}$ be a sequence of real independent and identically distributed random variables with $\mathbb{E}(X_{i})=0.$ Assuming also that there is a uniform bound for the $L_{2}$-norm and the supremum norm i.e. $\mathbb{E}(\rvert X_{i} \lvert)^2 \leq M a_{n} ,$ $\lvert \lvert X_{i}  \rvert \rvert_{\infty}  \leq M \leq n$ for every $i=1,...,n.$ Then for every $t\leq 2 \sqrt{M a_n}$ 
  \[ \mathbb{P}\big( \{ \frac{\lvert\sum_{i=1}^{n}X_{i}  \rvert   }{n} \geq t\} \big)\leq 2\exp{(-\frac{t^2 n}{M a_n})} .\]

  \begin{proof}

  \begin{align*}
         \mathbb{P}\biggl(\frac{1}{n} \sum_{i=1}^{n} X_{i} \geq t \biggr)  &= \mathbb{P}\biggl(\frac{\lambda}{n} \sum_{i=1}^{n} X_{i} \geq \lambda t \biggr) \\
    &= \mathbb{P}\biggl(\exp{(\frac{\lambda}{n} \sum_{i=1}^{n} X_{i})} \geq \exp{ (\lambda t )\biggr)   }\\
 &\leq  \exp{(-\lambda t)} \mathbb{E}\exp\biggl(\frac{\lambda}{n} \sum_{i=1}^{n} X_{i}  \biggl)\\
 &\leq  \exp{(-\lambda t)} \prod_{i=1}^{n} \mathbb{E}\exp\biggl(\frac{\lambda}{n}  X_{i}  \biggl).
  \end{align*}
  Where the inequalities are provided by  Chebysev inequality and the independence of the random variables. Moreover, for convenience, let $Y_j=\frac{Xj}{n}$ and we observe that, $ \lvert \lvert Y_{i}  \rvert \rvert_{\infty} \leq 1$ and $\mathbb{E} (Y_i)^2 \leq \frac{M a_n}{n^2 }=\sigma^2 $
\begin{align*}
     \mathbb{E}\exp\biggl(\frac{\lambda}{n}  X_{i}  \biggl) &= \mathbb{E} \biggl(1 +\sum_{k=2}^{\infty} \frac{\lambda^k Y_i^k}{ k!} \biggl)\\
     &= 1 +\sum_{k=2}^{\infty} \frac{\lambda^k \mathbb{E}( Y_i^k)}{ k!} \\
     &\leq 1 +\sum_{k=2}^{\infty} \frac{\lambda^k (\mathbb{E} Y_i^2)^\frac{k}{2} \lvert \lvert Y_{i}  \rvert \rvert_{\infty}^{\frac{k}{2}} }{ k!}\\
     &\leq 1 + \sum_{k=2}^{\infty} \frac{\lambda^k ((\sigma^2)^\frac{1}{2})^k}{k!}\\
     &=1+\exp{(\lambda\sigma)}-\lambda\sigma-1\\
     &\leq 1+\lambda \sigma+ (\lambda\sigma)^{2} -\lambda\sigma\\
     &\leq \exp{(\lambda\sigma)^{2}}
\end{align*}
where we have used that $\exp{(\lambda\sigma)} \leq  1+\lambda \sigma+ (\lambda\sigma)^{2} $ when $\lambda\sigma \leq 1$ and $1+x\leq e^x.$
Hence, 
\begin{align*}
\mathbb{P}\biggl(\frac{1}{n} \sum_{i=1}^{n} X_{i} \geq t \biggr) &\leq \exp{(-\lambda t)} \exp{(\lambda\sigma)^{2}n}\\
&\leq \exp{(-\frac{t^2}{2\sigma^{2} n})}\\
&=\exp{(-\frac{t^2 n}{M a_n})}
\end{align*}
where we choose $\lambda=\frac{t}{2\sigma^{2} n }$ which is an accepted value of $\lambda \iff t \leq 2 \sqrt{ M a_n}.$ We conclude the proof by replacing $X_i$ with $-X_i.$
  \end{proof}

  \end{lemma}
Now we move to the next necessary lemma. Here our target random variables are multiplied by centered independent Gaussian noises, but we can still obtain a similar tail bound by getting advantage the shape of the Gaussian tales.

\begin{lemma}\label{Bound_lemma2}

  Let $X_{1},..., X_{n}$ be a non-negative sequence of independent and identically distributed random variables with  $\mathbb{E}( X_{i} )^2 \leq M a_{n},$ $ \lvert \lvert X_{i}  \rvert  \rvert_{\infty} \leq M \leq n$ for every $i=1,...,n.$ Let also a sequence of independent random variables $\epsilon_{i}$ following normal distribution with zero mean and finite variance $\Tilde{\sigma}^2,$ for every $i=1,...,n. $ We assume also that $\epsilon_{i}$ are independent from $X_{i} $  for every $i=1,...,n. $\\
    Then for every $t\leq 2 \sqrt{M a_n}$ 
  \[ \mathbb{P}\big( \{ \frac{\lvert\sum_{i=1}^{n}X_{i}  \rvert   }{n} \geq t\} \big)\leq 2\exp{(-\frac{t^2 n}{M a_n \Tilde{\sigma}^2})} .\]

    \begin{proof}

First of all, we observe from the proof of Lemma 1 that, 
\[\mathbb{P}\biggl(\frac{1}{n} \sum_{i=1}^{n} X_{i} \epsilon_i \geq t \biggr) \leq   \exp{(-\lambda t)} \prod_{i=1}^{n} \mathbb{E}\exp\biggl(\frac{\lambda}{n}  X_{i}  \epsilon_i \biggl),\]
and

\begin{align*}
 \mathbb{E}\exp\biggl(\frac{\lambda}{n}  X_{i} \epsilon_i \biggl) &= \mathbb{E} \biggl(1 +\sum_{k=2}^{\infty} \frac{\lambda^k Y_i^k  \epsilon_i}{ k!} \biggl)\\
 &\leq 1 +\sum_{k=2}^{\infty} \frac{\lambda^k (\mathbb{E} Y_i^2)^\frac{k}{2} \lvert \lvert Y_{i}  \rvert \rvert_{\infty}^{\frac{k}{2}} \mathbb{E}(\epsilon_i)^k}{ k!}\\
 &\leq 1+ \sum_{k=2}^{\infty} \frac{\lambda^k ((\sigma^2)^\frac{1}{2})^k  \mathbb{E}(\epsilon_i)^k}{k!}\\
  &=  \mathbb{E}\exp\biggl(\lambda \sigma \epsilon_{i}\biggr)  \\
 &\leq    \mathbb{E}\exp\biggl(\lambda^2 \sigma^2  \Tilde{\sigma}^2\biggr)  
   \end{align*}

  Finally,
   \begin{align*}
\mathbb{P}\biggl(\frac{1}{n} \sum_{i=1}^{n} X_{i} \epsilon_i \geq t \biggr) &\leq \exp{(-\lambda t)} \exp{(\lambda\sigma \Tilde{\sigma})^{2}n}\\
&\leq \exp{(-\frac{t^2}{2\sigma^{2} \Tilde{\sigma}^2n})}\\
&=\exp{(-\frac{t^2 n}{M a_n})}
\end{align*}
where we choose $\lambda=\frac{t}{2\sigma^{2} \Tilde{\sigma}^2 n }$ which is an accepted value of $\lambda \iff t \leq 2 \sqrt{ M a_n}.$ We conclude the proof by replacing $X_i$ with $-X_i.$

    \end{proof}
\end{lemma}

\begin{theorem}

$   \textbf{Y}=m(\textbf{X}) +\epsilon$ where $\epsilon$ is a zero mean Gaussian noise with finite variance $\sigma $ independent of $\textbf{X} $. Assuming also that $\textbf{X} $ is uniformly distributed in $ [0,1]^d$ and $m$ is a Lipschitz function. Then there exists  constants $c_1, c_2$  depending on $d, \sigma $ and  $ \lVert m \rVert_{\infty} = \sup_{x \in [0,1]^d} | m(x) |$ such that,

\[  \mathbb{E}(\tilde{m}_{n,\infty}^{cc}(x) - m(x))^2 \leq
  c_{1} \left(1-\frac{1}{2d}\right)^{2\log_{2}{n}} 
  +C_3  \frac{ \log_2({ \log_2 n } )^d }{(\log_2 n)^{ \frac{d-1}{2}  }} \log{ 
 \bigg( \frac{ \log{n}  (\log_2 n)^{ \frac{d-1}{2}  }   }{  \log_2({ \log_2 n } )^d    }   \bigg)}    \]

    \begin{proof}
        Following the notation in \cite{S}, let $x \in [0,1]^d$, $ \lVert m \rVert_{\infty} = \sup_{x \in [0,1]^d} | m(x) |$,
and by the construction of the algorithm
\[ \tilde{m}_{n,\infty}^{Cen}(x) = \frac{\sum_{i=1}^{n}Y_{i}K_{k}(x,X_{i})}{\sum_{i=1}^{n}K_{k}(x,X_{i})} .\]
Let \[A_{n}(x) = \frac{1}{n}\sum_{i=1}^{n} \left( \frac{ Y_{i}K_{k}(x,X_{i}) - \mathbb{E}(Y K_{k}(x,X))} {\mathbb{E}(K_{k}(x,X))} \right),\]
\[B_{n}(x) = \frac{1}{n}\sum_{i=1}^{n} \left( \frac{K_{k}(x,X_{i}) - \mathbb{E}( K_{k}(x,X))} {\mathbb{E}(K_{k}(x,X))} \right),\]
and \[M_{n}(x) = \frac{\mathbb{E}(YK_{k}(x,X))}{\mathbb{E}( K_{k}(x,X))}.\]
Hence, we can reformulate the estimator as
\[ \tilde{m}_{n,\infty}^{Cen}(x) = \frac{M_{n}(x)+A_{n}(x)}{B_{n}(x)+1}. \]
Let $t \in (0,\frac{1}{2})$ and the event $C_{t}(x)$ where $\{A_{n}(x),B_{n}(x) \leq t\}.$
\begin{align*}
\mathbb{E}(\tilde{m}_{n,\infty}^{cc}(x) - m(x))^2 &= \mathbb{E}(\tilde{m}_{n,\infty}^{cc}(x) - m(x))^2 \mathbbm{1}_{C_{t}(x)} + \mathbb{E}(\tilde{m}_{n,\infty}^{cc}(x) - m(x))^2 \mathbbm{1}_{C_{t}^c(x)} \\
&\leq \mathbb{E}(\tilde{m}_{n,\infty}^{cc}(x) - m(x))^2 \mathbbm{1}_{C_{t}^c(x)} + c_{1} \left(1-\frac{1}{2d}\right)^{2k} + c_{2}t^2. \quad 
\end{align*}
Where the last inequality was obtained in \cite[p.1496]{S}
Moreover, in \cite{S},
\[ \mathbb{E}(\tilde{m}_{n,\infty}^{cc}(x) - m(x))^2 \mathbbm{1}_{C_{t}^c(x)} \leq c_{3}(\log{n})(\mathbb{P}(C_{t}^c(x)))^{\frac{1}{2}}.\]

\begin{proposition}
    
The probability $ \mathbb{P}(C_{t}^c(x)) \leq C \exp{( -\frac{t^2 n}{2^k a_n} )}$ for a constant $C$ independent of $k,n$ and $a_n$ is a sequence that converges to zero as $n$ tends to infinity.

    \begin{proof}
First of all, we notice that  \[  \mathbb{P}(C_{t}^c(x)) \leq    \mathbb{P}( \rvert A_{n}(x) \rvert >t)+  \mathbb{P}( \lvert B_{n}(x) \rvert >t) .\]

 The related result will follow by working separately on both inequalities.

By lemma B.4 \cite{Interp} for all $x \in [0,1]^d $ and $ d\geq 2$ we have \[  \mathbb{E}\big( K^c_k(x,X) \big)^2 \leq  \frac{C_1+C_2(\log_2(k))^d}{k^{\frac{d-1}{2} }2^k} \]  where
$C_1=1+ \frac{2d^{\frac{d}{2}}}{ (4 \pi)^{\frac{d-1}{2}} }$ , $ C_2=5^d(\frac{d-1}{2})^d$ and for convenience let $a_n= \frac{C_1+C_2(\log_2(k))^d}{k^{\frac{d-1}{2} }}.$

Hence,  let $  \tilde{X}_{i} =   \frac{K_{k}(x,X_{i})}{\mathbb{E}(K_{k}(x,X))} -1$ a sequence of centered iid random variables with \[\lvert \lvert  \tilde{X}_{i}  \rvert \rvert_{\infty}=\sup \{  \lvert \tilde{X}_{i} \rvert\}= \sup \{  \lvert \frac{K_{k}(x,X_{i})}{\mathbb{E}(K_{k}(x,X))} -1 \rvert\} \leq \frac{1}{\mathbb{E}(K_{k}(x,X))} \sup{ K_{k}(x,X_{i})}+1  \leq  2^{k+1}\] and

\begin{align*}
     \mathbb{E}\big( \frac{K^c_k(x,X_i )}{ \mathbb{E} (K^c_k(x,X )} \big)^2 &= \frac{1}{ (\mathbb{E} (K^c_k(x,X )))^2}\mathbb{E}\big( K^c_k(x,X) \big)^2\\
     &\leq \frac{1}{(\frac{1}{2^k})^2}  \frac{C_1+C_2(\log_2(k))^d}{k^{\frac{d-1}{2} }2^k}\\
     &= 2^k a_n
\end{align*}
By lemma \ref{Bound_lemma1}, for every $t \leq 2\sqrt{2^k a_n}$
\[ \mathbb{P}( \lvert B_{n}(x) \rvert >t) \leq 2\exp{(-C\frac{  t^2 n}{2^k  a_n})}  \]
We need an estimate for the  $ \mathbb{P}( \rvert A_{n}(x) \rvert >t)$ where, \[A_{n}(x)= \frac{1}{n}\sum_{i=1}^{n} \big( \frac{ Y_{i}K_{k}(x,X_{i})  - \mathbb{E}(Y K_{k}(x,X))} {\mathbb{E}(K_{k}(x,X))}  \big)   .\]

With simple calculations,    
\begin{align*}
    A_{n}(x) &= \frac{1}{n}\sum_{i=1}^{n} \bigg( \frac{ Y_{i}K_{k}(x,X_{i})  - \mathbb{E}(Y K_{k}(x,X))} {\mathbb{E}(K_{k}(x,X))}  \bigg)\\
    &=\frac{1}{n}\sum_{i=1}^{n} \bigg( \frac{ m(X_{i})K_{k}(x,X_{i})  - \mathbb{E}(m(X) K_{k}(x,X))} {\mathbb{E}(K_{k}(x,X))}  \bigg) +\frac{1}{n}\sum_{i=1}^{n} \bigg( \frac{ \epsilon_{i}K_{k}(x,X_{i})  - \mathbb{E}(\epsilon K_{k}(x,X))} {\mathbb{E}(K_{k}(x,X))}  \bigg)\\
    &=\frac{1}{n}\sum_{i=1}^{n} \bigg( \frac{ m(X_{i})K_{k}(x,X_{i})  - \mathbb{E}(m(X) K_{k}(x,X))} {\mathbb{E}(K_{k}(x,X))}  \bigg) +\frac{1}{n}\sum_{i=1}^{n} \bigg( \frac{ \epsilon_{i}K_{k}(x,X_{i}) }{\mathbb{E}(K_{k}(x,X))} \bigg) .
\end{align*}

Let $Z_i= 2 \frac{ \epsilon_{i}K_{k}(x,X_{i}) }{\mathbb{E}(K_{k}(x,X))} $ a sequence of centered iid random variables with

\begin{align*}
    \mathbb{E} (Z_i)^2 
     &\leq \frac{4}{\frac{1}{(2^k)^2}}   \mathbb{E} (\epsilon K_{k}(x,X))^2 \\
     &=  \frac{4}{\frac{1}{(2^k)^2}} \mathbb{E}(\epsilon)^2 \mathbb{E}(K_{k}(x,X)))^2) \\
     &\leq \tilde{c}  2^k  \sigma^2  a_n
\end{align*}

and \[ \lvert \lvert    \frac{ K_{k}(x,X_{i}) }{\mathbb{E}(K_{k}(x,X))}  \rvert \rvert_{\infty} \leq     2^k \]

hence by \ref{Bound_lemma2} \[ \mathbb{P} \big( \lvert Z_i \rvert \geq t \big)=\mathbb{P} \big( \frac{2}{n}\sum_{i=1}^{n} \lvert \frac{ \epsilon_{i}K_{k}(x,X_{i}) } {\mathbb{E}(K_{k}(x,X))} \rvert \geq t  \big) \leq 2\exp{(-\frac{t^2 n}{2^k  a_n}}) \]  for every $t \leq 2\sqrt{2^k a_n}$

\begin{align*}
    \mathbb{P} \big( \lvert A_{n}(x) \rvert \geq t \big)  &\leq 
    \mathbb{P} \bigg( \bigg\lvert \frac{2}{n} \sum_{i=1}^{n} \frac{m(X_{i})K_{k}(x,X_{i}) - \mathbb{E}(m(X) K_{k}(x,X))}{\mathbb{E}(K_{k}(x,X))} \bigg\rvert \geq t \bigg) \\
    &\quad+ \mathbb{P} \bigg( \bigg\lvert \frac{2}{n} \sum_{i=1}^{n} \frac{\epsilon_{i}K_{k}(x,X_{i})}{\mathbb{E}(K_{k}(x,X))} \bigg\rvert \geq t \bigg)\\
    &\leq  2 \exp{(-C_1\frac{t^2 n}{2^k  }}) + 2\exp{(-C_2\frac{t^2 n}{2^k  a_n}}) 
    \leq C \exp{(-C_3\frac{t^2 n}{2^k  a_n}}) .
\end{align*}
where we have used  a partial result from \cite{I-A} [proposition 6.]
and finally,

\[\mathbb{P}(C_{t}^c(x)) \leq 2 \exp{(-\tilde{C} \frac{t^2 n}{2^k a_n} )} .\]
    \end{proof}
\end{proposition}

To obtain the desired rate of convergence, we need an upper bound for
\[\mathbb{E}(\tilde{m}_{n,\infty}^{cc}(x) - m(x))^2 \leq  c_3 \log{n} \exp{(-\tilde{C} \frac{t^2 n}{2^k a_n} )} + c_{1} \left(1-\frac{1}{2d}\right)^{2k} + c_{2}t^2  \]
and we choose $2^k=n$ in the mean interpolation regime, 
\[ \mathbb{E}(\tilde{m}_{n,\infty}^{cc}(x) - m(x))^2 \leq  c_3 \log{n} \exp{(-\tilde{C} \frac{t^2 }{ a_n} )} + c_{1} \left(1-\frac{1}{2d}\right)^{2\log_{2}{n}} + c_{2}t^2\]

Finally, by minimizing the right hand of the equation in terms of $t$ one has that, 
$t^2= Ca_n \log{\frac{\log{n}}{a_n}} $ 
and of course,
\begin{align*}
     c_3 \log{n} \exp{(-\tilde{C} \frac{t^2 }{ a_n} )} + c_{1} \left(1-\frac{1}{2d}\right)^{2\log_{2}{n}} + c_{2}t^2 &\leq c_{1} \left(1-\frac{1}{2d}\right)^{2\log_{2}{n}}\\
     &\quad  + C_3 a_n + C_2 a_n \log{(  \frac{\log{n}}{a_n}  )}
\end{align*}

and therefore, since $ a_n = \frac{ \log_2({ \log_2 n } )^d }{(\log_2 n)^{ \frac{d-1}{2}  }     }$

\[  \mathbb{E}(\tilde{m}_{n,\infty}^{cc}(x) - m(x))^2 \leq
  c_{1} \left(1-\frac{1}{2d}\right)^{2\log_{2}{n}} 
  +C_3  \frac{ \log_2({ \log_2 n } )^d }{(\log_2 n)^{ \frac{d-1}{2}  }} \log{ 
 \bigg( \frac{ \log{n}  (\log_2 n)^{ \frac{d-1}{2}  }   }{  \log_2({ \log_2 n } )^d    }   \bigg)}    \]
 and this concludes the proof.
    \end{proof}
    
\end{theorem}

\bibliography{main.bib}

\bibliographystyle{acm}

\end{document}